\documentclass[german,a4paper,12pt]{article}
\usepackage{amsmath,amsthm,amsfonts,latexsym,amscd}

\usepackage{amsthm}

\usepackage{amssymb}
\swapnumbers
\theoremstyle{plain}

\theoremstyle{definition}

\theoremstyle{remark}

\begin{document}

\centerline{\bf Semi-local Zariski-van Kampen theorems}
\vskip.2in
\centerline{\bf Helmut A. Hamm}
\vskip.2in

{\bf Abstract:} The classical second Lefschetz theorem is a complement to the well-known Lefschetz theorem on hyperplane sections for smooth complex projective varieties. When passing from projective varieties to quasi-projective ones it is in view of the work of Zariski and van Kampen, who started with the case of the complement of a projective hypersurface, preferable to speak of a Zariski-van Kampen theorem instead, as done in the lecture notes by C. Eyral [E]. The present paper starts from some local analogue developed in [H2] - both for homology and homotopy groups - and passes from the local to the semi-local case. The reason is that as a consequence we obtain a deeper understanding of the local case, by including the nearby fibre, as well as a bridge to the global case (of quasi-projective varieties).\\

First let us recall classical achievements of Lefschetz and Zariski, starting with the Lefschetz theorem on hyperplane sections: Suppose that $X\subset \mathbb{P}^N(\mathbb{C})$ is a smooth subvariety of pure dimension $n>0$ and $L$ a hyperplane in $\mathbb{P}^N$ (which is general enough but this is here not necessary). Then $H_k(X,X\cap L)=0$ for $k<n$, i.e. $H_k(X\cap L)\to H_k(X)$ is bijective for $k<n-1$ and surjective for $k=n-1$ - with homology we always mean integral homology. See [L]. In his proof, Lefschetz used a pencil method.\\
In fact, one can say more, as found out later on: the pair $(X,X\cap L)$ is $(n-1)$-connected, so $\pi_k(X\cap L,\cdot)\to \pi_k(X,\cdot)$ is bijective for $k<n-1$ and surjective for $k=n-1$, $\cdot$ being any base point in $X\cap L$.\\
It is possible to replace $X$ by a smooth quasi-projective variety $X=\overline{X}\setminus Z$, $Z\subset\overline{X}$ being projective varieties, see [HL1]. But then $L$ must be supposed to be chosen general. In a special case this goes back to Zariski [Z], so we speak of a Lefschetz-Zariski theorem. This works for homology and homotopy groups.\\
As said above $H_{n-1}(X\cap L)\to H_{n-1}(X)$ is merely surjective, so one has to determine the kernel in order to compute $H_{n-1}(X)$ from $H_{n-1}(X\cap L)$. If $X$ is projective this has been done by Lefschetz using his pencil method: the kernel can be computed with the help of some variation mapping. This works similarly for homotopy groups.\\
Zariski and van Kampen, who were interested in the fundamental group, transferred the pencil method to a special type of a quasi-projective variety where we deal with the complement of some hypersurface. Therefore theorems about the kernel above are called Zariski-van Kampen theorems.\\
A general homological Zariski-van Kampen theorem for smooth quasi-projective varieties has been proved by Chéniot [C]. A corresponding theorem for homotopy groups has been conjectured by Chéniot and Eyral [CE].\\
In (local) singularity theory, when passing from hypersurfaces to complete intersections the author found decades ago a local Lefschetz theorem on hyperplane sections [H1]. This was the starting point for local theorems analogous to the global ones before [HL2]. Now in [H2] a local Zariski-van Kampen theorem has been proved both in the framework of homology and homotopy groups. In fact the statement in the case of homotopy groups is inspired by the conjecture of Chéniot and Eyral mentioned before.\\
How to pass from the local theorem to a global one for quasi-projective varieties, especially in the case of homotopy groups? In this paper we propose a semi-local generalization as a bridge: it includes the local case as well as a preparation of the global one.\\
Main results are a semi-local Zariski-van Kampen theorem for  homology (Theorem 1.2.3) and homotopy groups (Theorem 2.2.2).\\

\vskip.2in
\centerline{\bf Contents}
\vskip.1in

0.   Main results\\

1.   Homology\\
1.0. A local Zariski-van Kampen theorem\\
1.1. Towards the semi-local case\\
1.2. The semi-local case\\
1.3. First application: hyperplane section of the local nearby fibre\\
1.4. Consequences for local homology groups\\
1.5. Second application: the global case\\
1.6. Additional comments\\

2.   Homotopy groups\\
2.0. Preparation\\
2.1. Towards the semi-local case\\
2.2. The semi-local case\\
2.3. Application: hyperplane section of the local nearby fibre\\
2.4. Consequences for local homotopy groups\\

\vskip.2in

\centerline{\bf 1. Homology}
\vskip.1in

{\bf 1.0. A (strong) local Zariski-van Kampen theorem}\\

In [H2] we proved a strong local Zariski-van Kampen theorem for homology - as said before we always mean integral homology. Before recalling it let us make the data precise:\\

Let $U$ be an open neighbourhood of $0$ in $\mathbb{C}^N$, $Z$ and $\bar{X}$ closed complex analytic subsets of $U$, $0\in Z\subset \bar{X}$, $X:=\bar{X}\setminus Z$.\\
Choose a Whitney stratification of $U$ which is compatible with $\bar{X}$ and $Z$.\\

A natural condition will be:\\

$(C_n)$: $X$ is smooth, $\dim(X,x)\ge n$ for all $x\in X$.\\

Note that $(C_n)$ implies the condition\\

$(C'_n)$: $rHd(X)\ge n$, where $rHd$ means the rectified homological depth, cf. [HL3] 1.1.3. This means that\\ 
for every $m\ge -1$ there is a complex analytic subset $X_m$ of $X$ of dimension $\le m$ such that $H_{n+k}(X,X\setminus\{x\})=0$ for all $x\in X\setminus X_m, k\le m$,\\ 
i.e. (using K\"unneth) that\\ 
$H_{n-j-1}(V,V\setminus S)=0$, $j\le k$, if $S$ is a stratum of $X$ of dimension $k$, $x\in S$ and $V$ belongs to a fundamental system of neighbourhoods of $x$ in $X$.\\
For instance, this condition holds if $X$ is a local complete intersection which is everywhere of dimension $\ge n$.\\

If $n\ge 1$, $(C_n)$ implies the condition\\

$(C''_n)$: As in the case of {\bf $(C'_n)$}, with $X_0=\emptyset$.\\

Also, $(C'_n)$ implies $(C''_n)$ if all strata of $X$ have dimension $>0$.\\

In subsection 1.0 we suppose that $(C_n)$ or at least $(C'_n)$ holds.\\

Let $0<\varepsilon\ll 1$. Put ${\rm B}_\varepsilon:=\{z\in\mathbb{C}^N\,|\,\|z\|\le\varepsilon\}$. Then $\partial {\rm B}_\varepsilon$ intersects the strata transversally. Let $f:U\to\mathbb{C}$ be a holomorphic function with an isolated singularity at $0$, $f(0)=0$. Then $f^{-1}(0)\setminus\{0\}$ is transverse to $\partial {\rm B}_\varepsilon\cap\bar{X}$ with respect to the induced stratification. Choose $0<\alpha\ll\varepsilon$. Put\\
$X_\varepsilon:=X\cap{\rm B}_\varepsilon,  \partial X_\varepsilon:=X\cap\partial{\rm B}_\varepsilon$, $F_t:=f^{-1}(t)\cap X_\varepsilon$, $\partial F_t:=f^{-1}(t)\cap \partial X_\varepsilon$,\\
and $F:=F_\alpha$, $\partial F:=\partial F_\alpha$.\\
Now recall:\\

{\bf Theorem 1.0.1 (Strong local Lefschetz-Zariski theorem} (cf. [HL2], [HL3]): $H_k(X_\varepsilon,F)=0, k<n$,
so $H_k(F)\stackrel{\simeq}{\to} H_k(X_\varepsilon), k\le n-2$,\\  
and $H_{n-1}(F)\to H_{n-1}(X_\varepsilon)$ is surjective.\\

Now note that $X_\varepsilon$ admits the deformation retract $\partial (X_\varepsilon\cap\{|f|\le\alpha\})=(X_\varepsilon\cap\{|f|=\alpha\})\cup((\partial X_\varepsilon)\cap\{|f|\le\alpha\})$ - recall that $0\notin X$. As shown in [H2], we have that $g=\frac{f}{\alpha}:\partial (X_\varepsilon\cap\{|f|\le\alpha\})\to \mathbb{D}_2$ gives a stratified open book as defined in [H2], where $\mathbb{D}_2:=\{t\in\mathbb{C}\,|\,|t|\le 1\}$:\\
In fact, let $M$ be a suitable open neighborhood of $\partial(\bar{X}_\varepsilon\cap\{|f|\le\alpha\})$ in $\mathbb{C}^N$. Choose a stratification of $M$ such that\\
$A:=\bar{X}_\varepsilon\cap M, B:=Z\cap A, A':=\partial\bar{X}_\varepsilon\cap M$ and the intersections of these sets are unions of strata. Then we may achieve the conditions for a stratified open book, with\\
$X^1:=X_\varepsilon\cap\{|f|=\alpha\}$, $X^2:=\partial X_\varepsilon\cap\{|f|\le\alpha\}$, $g:=\frac{f}{\alpha}$. In particular, there is a corresponding geometric monodromy $h:F\to F$, which implies that $h$ is a monodromy for $\frac{f}{\alpha}:X^1\to S^1$ with $h|\partial F=id$, see [H2] section 5b). So we have a corresponding variation\\
$var_k(\frac{f}{\alpha}):H_k(F,\partial F)\to H_k(F)$ defined by\\
$var_k(\frac{f}{\alpha})([c]):=[h\circ c-c]$, see [H2] section 4.\\

On the other hand, we have, using [H2] Theorem 6.2 (without the assumption $(C_n)$ or $(C'_n)$):\\

{\bf Theorem 1.0.2:} a) $w_k:H_k(F,\partial F)\to H_{k+1}(X_\varepsilon,F)$ is an isomorphism, where $C:=\partial X_\varepsilon\cap\{|f|=\alpha,Im\,f\ge 0\}$ and $w_k$ is the composition of the isomorphisms
$$\begin{array}{c}
H_k(F,\partial F)\\
\downarrow\\
H_k((X_\varepsilon\cap\{f=\pm\alpha\})\cup C,(X_\varepsilon\cap \{f=-\alpha\})\cup C)\\
\uparrow\\
H_{k+1}(X_\varepsilon\cap\{|f|=\alpha,Im\,f\ge 0\}, (X_\varepsilon\cap\{f=\pm\alpha\})\cup C)\\
\downarrow\\ 
H_{k+1}(X_\varepsilon,(X_\varepsilon\cap\{|f|=\alpha, Im\,f\le 0\})\cup (\partial X_\varepsilon\cap\{|f|\le\alpha\}))\\
\uparrow\\
H_{k+1}(X_\varepsilon,F)
\end{array}$$
(Warning: $w_k$ has nothing to do with the exact homology sequence of the triple $(X_\varepsilon,F,\partial F)$ !)\\
b) $\partial\circ w_k:H_k(F,\partial F)\to H_k(F)$ coincides with $-var_k(\frac{f}{\alpha})$.\\

As a consequence of Theorem 1.0.1 and 1.0.2 we get:\\

{\bf Corollary 1.0.3 (Local strong Zariski-van Kampen theorem)} (see [H2] Corollary 6.3): 
$$H_{n-1}(X_\varepsilon)\simeq H_{n-1}(F)/Im\,var_{n-1}(\frac{f}{\alpha})$$\\

{\bf 1.1. Towards the semi-local case}\\

Now we would like to pass from the local to a semi-local situation.  
First we consider a simpler case: semi-local with respect to the source and local with respect to the target space. This means: drop the condition of section 1.0 that $\varepsilon$ is small, so replace $\varepsilon$ by $r>0$, we will use $\varepsilon$ for something else.\\
So let $U$ be an open subset of $\mathbb{C}^N$, $r>0$, ${\rm B}_r:=\{z\in\mathbb{C}^N\,|\,\|z\|\le r\}\subset U$, $f:U\to \mathbb{C}$ holomorphic. Let $\bar{X}$ and $Z$ be complex analytic subsets of $U$, $X=\bar{X}\setminus Z$. Choose a complex analytic Whitney stratification of $U$ which is compatible with $\bar{X}$ and $Z$. We assume that $\{\|z\|=r\}$ is transverse to this stratification, so we obtain a refined stratification such that $\bar{X}\cap\partial{\rm B}_r$ is a union of strata.\\ 
We assume furthermore that $f^{-1}(0)$ is a space of dimension $N-1$, smooth along $\bar{X}\cap\partial{\rm B}_r$, which intersects the strata of this space transversally. \hfill{(*)}\\

So there are finitely many points $z_1,\ldots,z_k\in\bar{X}\cap{\rm B}_r\cap\{f=0\}$ for which $\|z_j\|<r$, $f^{-1}(0)$ not transversal to $(\bar{X},Z)$ at $z_j$. Let $0<\beta\ll \varepsilon\ll r$, $X_r:=X\cap {\rm B}_r$, $F_\beta:=f^{-1}(\beta)\cap X_r$. Put $B_j:=\{z\in\mathbb{C}^N\,|\,\|z-z_j\|\le\varepsilon\}$.\\

{\bf Lemma 1.1.1:} For every $m$ we have an isomorphism
$$\oplus_{j=1}^kH_m(X_r\cap B_j\cap\{|f|\le\beta\},F_\beta\cap B_j)\stackrel{\simeq}{\to} H_m(X_r\cap\{|f|\le\beta\},F_\beta)$$

{\bf Proof:} On the right, we may replace $X_r\cap\{|f|\le\beta\}$ by its deformation retract $F_\beta\cup((B_1\cup\ldots\cup B_k)\cap X\cap\{|f|\le\beta\})$. Then apply excision.\\
More precisely, with $B:=B_1\cup\ldots\cup B_k$:\\
$F_\beta\cap\partial B$ is a deformation retract of $X\cap\{|f|\le\beta\}\cap\partial B$\\
(we have a fibration over $\{|t|\le\beta\}$),\\
hence $F_\beta\setminus\stackrel{o}{B}$ is a deformation retract of
$(F_\beta\setminus\stackrel{o}{B})\cup (X\cap\{|f|\le\beta\}\cap\partial B)$,\\
and simultaneously of $X_r\cap\{|f|\le\beta\}\setminus\stackrel{o}{B}$\\
(which is the total space of a fibration over $\{|t|\le\beta\}$),\\
so $(F_\beta\setminus\stackrel{o}{B})\cup (X\cap\{|f|\le\beta\}\cap\partial B)$ is a deformation retract of $X\cap\{|f|\le\beta\}\setminus\stackrel{o}{B}$.\\
Hence $F_\beta\cup(X\cap\{|f|\le\beta\}\cap B)$ is a deformation retract of $X_r\cap\{|f|\le\beta\}$, by excision.\\
  
As for the terms on the left, we have for any $m$:\\

{\bf Lemma 1.1.2:} a) If $z_j\in X$, we have an isomorphism
$$H_m(X_r\cap B_j\cap\{|f|\le\beta\},F_\beta\cap B_j)\stackrel{\partial}{\to} \widetilde{H}_{m-1}(F_\beta\cap B_j)$$
b) If $z_j\not\in X$,
$$H_m(X_r\cap B_j\cap\{|f|\le\beta\},F_\beta\cap B_j)\stackrel{\simeq}{\to} H_{m-1}(F_\beta\cap B_j,F_\beta\cap\partial B_j)$$
the mapping being explained below.\\

{\bf Proof:} a) Look at the long exact homology sequence of the triad 
$(X_r\cap B_j\cap\{|f|\le\beta\},F_\beta\cap B_j,*)$, $*$ being a point of $F_\beta\cap B_j$, and note that $*$ is a deformation retract of the contractible space $X_r\cap B_j\cap\{|f|\le\beta\}$.\\
b) See Theorem 1.0.2a), the inverse of the mapping being called $w_{m-1}$ there.\\

Remember that we have introduced conditions $(C_n)$, $(C'_n)$ and $(C''_n)$ in subsection 1.0.\\

{\bf Lemma 1.1.3:} Suppose that we have condition $(C_n)$ or at least $(C'_n)$. Then   $H_m(X_r\cap\{|f|\le\beta\},F_\beta)=0, m\le n-1$.\\

{\bf Proof:} Use Lemma 1.1.1 and the strong local Lefschetz-Zariski theorem (Theorem 1.0.1).\\

{\bf Theorem 1.1.4:} We assume $(C_n)$ or at least $(C'_n)$. Then
$$H_{n-1}(F_\beta)\to H_{n-1}(X_r\cap\{|f|\le\beta\})$$
is an epimorphism whose kernel 
is generated by $N_1,\ldots,N_k$, where $N_j$ is the image of the mapping:\\
a) $\widetilde{H}_{n-1}(F_\beta\cap B_j)\to H_{n-1}(F_\beta\cap B_j)\to H_{n-1}(F_\beta)$ if $z_j\in X$ resp.\\
b) $H_{n-1}(F_\beta\cap B_j,F_\beta\cap\partial B_j)\to H_{n-1}(F_\beta\cap B_j)\to H_{n-1}(F_\beta)$ if $z_j\not\in X$,\\ 
the first mapping being the variation (see Theorem 1.0.2b):\\$[c]\mapsto  [h_j\circ c-c]$,\\
with $h_j$ being a geometric monodromy of $\frac{f}{\beta}|(F_\beta\cap B_j)\cup(\partial B_j\cap X\cap\{|f|\le\beta\})$.\\

({\bf Remark:} For given $j$, we may alternatively use b) instead of a) if $z_j$ is contained in a stratum of $X$ of positive dimension, replacing $X$ by $X\setminus\{z_j\}$).\\

{\bf Proof:} By Lemma 1.1.3 we have an epimorphism indeed, and the kernel $N$ is the image of $H_n(X_r\cap\{|f|\le\beta\},F_\beta)\to H_{n-1}(F_\beta)$. \\
By Lemma 1.1.1, $N=\widetilde{N}_1+\ldots+\widetilde{N}_k$, where $\widetilde{N}_j$ is the image of\\ $H_n(X\cap B_j\cap\{|f|\le\beta\},F_\beta\cap B_j)\to H_n(X_r\cap\{|f|\le\beta\}),F_\beta)\stackrel{\partial}{\to} H_{n-1}(F_\beta)$, \\
or of $H_n(X\cap B_j\cap\{|f|\le\beta\},F_\beta\cap B_j)\stackrel{\partial}{\to} H_{n-1}(F_\beta\cap B_j)\to H_{n-1}(F_\beta)$:\\
use the commutative diagram
$$\begin{array}{ccc}
H_n(X\cap B_j\cap\{|f|\le\beta\},F_\beta\cap B_j)&\to&H_n(X_r\cap\{|f|\le\beta\}),F_\beta)\\
\downarrow&&\downarrow\\
H_{n-1}(F_\beta\cap B_j)&\to&H_{n-1}(F_\beta)
\end{array}$$
a) In this case we have a commutative diagram
$$\begin{array}{ccc}
\widetilde{H}_{n-1}(F_\beta\cap B_j)&&\\
\uparrow&\searrow&\\
H_n(X_r\cap B_j\cap\{|f|\le\beta\},F_\beta\cap B_j)&\to& H_{n-1}(F_\beta\cap B_j)
\end{array}$$
where the vertical arrow is bijective, cf. Lemma 1.1.2a).
So $\widetilde{N}_j=N_j$.\\
b) Similar, with the commutative diagram$$\begin{array}{ccc}
H_{n-1}(F_\beta\cap B_j, F_\beta\cap\partial B_j)&&\\
\uparrow&\searrow&\\
H_n(X_r\cap B_j\cap\{|f|\le\beta\},F_\beta\cap B_j)&\to& H_{n-1}(F_\beta\cap B_j)
\end{array}$$
where the first vertical is bijective, cf. Lemma 1.1.2b), and the diagonal is given by the variation, up to sign, cf. Theorem 1.0.2b). So $\widetilde{N}_j=N_j$ again.\\

Now we would like to avoid looking to two cases a) and b). It was observed already by Tib\v{a}r (cf. [Ti] Theorem 10.2.3) that condition $(C'_n)$ is then no longer sufficient and that we need a further condition on the local homology.\\ 

{\bf Theorem 1.1.5:} Suppose that we have condition $(C_n), n\ge 1$, or at least $(C''_n)$. Then the kernel discussed in Theorem 1.1.4 is the image of the variation $H_{n-1}(F_\beta,\partial F_\beta)\to H_{n-1}(F_\beta)$, defined by $[c]\mapsto [h(c)-c]$.\\

{\bf Proof:} We may assume that $z_1,\ldots,z_k$ are contained in $Z$: otherwise put $Z':=Z\setminus\{z_1,\ldots,z_k\}$, $X'_r:=X_r\setminus\{z_1,\ldots,z_k\}$. And\\
$H_{n-1}(X'_r\cap\{|f|\le\beta\})\to H_{n-1}(X_r\cap\{|f|\le\beta\})$ is bijective, since \\$H_m(X_r\cap\{|f|\le\beta\},X'_r\cap\{|f|\le\beta\})=0$ for $m=n-1,n$ because of condition $C''_n$.\\
So the kernel is the image of $H_{n-1}(F_\beta,F_\beta\setminus\bigcup_{j=1}^k \stackrel{o}{B}_j)\to H_{n-1}(F_\beta)$, since the left hand group is isomorphic to $\oplus H_n(F_\beta\cap B_j,F_\beta\cap\partial B_j)$. And\\
$F_\beta\setminus\bigcup_{j=1}^k \stackrel{o}{B}_j$ is homeomorphic to $F_0\setminus\bigcup_{j=1}^k \stackrel{o}{B}_j$,\\
where $F_0\setminus\bigcup_{j=1}^k \stackrel{o}{B}_j)$ is a deformation retract of $F_0\setminus\{z_1,\ldots,z_k\}$. Then apply Morse theory to mapping $z\mapsto -\|z\|^2$ in order to show that the pair\\
$(F_\beta\setminus\{z_1,\ldots,z_k\},\partial F_\beta)$ is $(n-2)$-connected, which implies that\\
$H_{n-1}(F_\beta,\partial F_\beta)\to H_{n-1}(F_\beta,F_\beta\setminus\bigcup_{j=1}^k \stackrel{o}{B}_j)$ is surjective - use the exact homology sequence of the triple\\
$(F_\beta,F_\beta\setminus\bigcup_{j=1}^k \stackrel{o}{B}_j,\partial F_\beta)$.\\

{\bf 1.2. The semi-local case}\\

We start from the situation studied in 1.1 but we replace the condition (*) by the following one, where $\alpha>0$:\\
If $|t|\le\alpha$, $f^{-1}(t)$ is smooth of dimension $N-1$ up to isolated singularities which do not lie on $\partial{\rm B}_r$, and it intersects the strata of $X\cap\partial{\rm B}_r$ transversally.\hfill{(**)}\\

Then there are finitely many points $z_1,\ldots,z_s$ in $X_r\cap\{|f|\le\alpha\}$ such that $f^{-1}(f(z_j))$ does not intersect $\bar{X}$ transversally at $z_j$ in the stratified sense. We assume that $|f(z_j)|<\alpha$ for all $j$. Note that $f(z_1),\ldots,f(z_s)$ are not nessarily distinct, let\\
$\{t_1,\ldots,t_k\}:=\{f(z_1),\ldots,f(z_s)\}$ if the set on the right has $k$ elements.\\
Let $0<\beta\ll\varepsilon\ll\min(\alpha,r)$ and $B_j:=\{z\in \mathbb{C}^N\,|\,\|z-z_j\|\le\varepsilon\}$.\\
For $j=1,\ldots,k$ let $D_j:=\{t\in\mathbb{C}\,|\,|t-t_j|\le\beta\}$, $\Gamma_j$ the trace of a simple path from $\alpha$ to $t_j+\beta$ which lies in\\
$\{t\in\mathbb{C}\,|\,|t|<\alpha, |t-t_1|>\beta,\ldots,|t-t_k|>\beta\}$ except for the initial and end point such that $\Gamma_1,\ldots,\Gamma_k$ do not intersect except at the initial point, $\Gamma:=\Gamma_1\cup\ldots\cup\Gamma_k$.\\
Let $\Delta_j:=D_j\cup \Gamma_j$, $\Delta:=\Gamma\cup D_1\cup\ldots\cup D_k$.\\
Note that we have induced isomorphisms
$$H_m(F_\alpha)\stackrel{\simeq}{\to} H_m(f^{-1}(\Gamma_j)\cap X)\stackrel{\simeq}{\leftarrow} H_m(F_{t_j+\beta})$$ 
because $f|f^{-1}(\Gamma_j)\cap X\to \Gamma_j$ is a trivial fibration.\\

{\bf Lemma 1.2.1:} For every $m$, we have an isomorphism
$$\oplus_{j=1}^k H_m(X_r\cap \{|f-t_j|\le\beta\},F_{t_j+\beta})\to H_m(X_r\cap\{|f|\le\alpha\},F_\alpha)$$

{\bf Proof:} It is obtained by composing the following isomorphisms:\\
$\oplus_{j=1}^k H_m(X_r\cap \{|f-t_j|\le\beta\},F_{t_j+\beta})\stackrel{\simeq}{\to} H_m(X_r\cap f^{-1}(\Delta),X_r\cap f^{-1}(\Gamma))$, using excision, and 
$$\begin{array}{ccc}
H_m(X_r\cap f^{-1}(\Delta),X_r\cap f^{-1}(\Gamma))&\stackrel{\simeq}{\leftarrow}&H_m(X_r\cap f^{-1}(\Delta),F_\alpha)\\
&&\downarrow\simeq\\
&&H_m(X_r\cap\{|f|\le\alpha\},F_\alpha)
\end{array}$$
because $X_r\cap f^{-1}(\Delta)$ is a deformation retract of\\
$X_r\cap\{|f|\le\alpha\}$ and $F_\alpha$ is a deformation retract of $X_r\cap f^{-1}(\Gamma)$.\\
In fact: $\Delta\setminus \stackrel{o}{D}_1\cup\ldots\cup\stackrel{o}{D}_k$ is a deformation retract of $\{t\,|\,|t|\le\alpha, t\notin \stackrel{o}{D}_1,\ldots,t\notin\stackrel{o}{D}_k\}$.\\
The same holds for the inverse images with respect to $f|X_r$. Therefore
$X_r\cap f^{-1}(\Delta)$ is a deformation retract of $X_r\cap\{|f|\le\alpha\}$.\\

Recall the conditions $C_n, C'_n, C''_n$ of subsection 1.0. \\

{\bf Lemma 1.2.2:} Suppose that we have $(C_n)$ or at least $(C'_n)$. Then $H_m(X_r\cap\{|f|\le \alpha\},F_\alpha)=0, m\le n-1$.\\

{\bf Proof:} This follows from Lemma 1.2.1 and Lemma 1.1.3.\\

{\bf Theorem 1.2.3 (Semi-local Zariski-van Kampen theorem for homology):} Suppose that we have $(C_n)$ or at least $(C'_n)$. Then the mapping
$$H_{n-1}(F_\alpha)\to H_{n-1}(X_r\cap\{|f|\le\alpha\})$$ is an epimorphism whose kernel is generated by the spaces $N_j$, where $N_j$ is the image of the mapping:\\
a) $\widetilde{H}_{n-1}(F_{f(z_j)+\beta}\cap B_j)\to H_{n-1}(F_{f(z_j)+\beta}\cap B_j)\to H_{n-1}(F_{f(z_j)+\beta})\to H_{n-1}(F_\alpha)$ if $z_j\in X$ resp.\\
b) $H_{n-1}(F_{f(z_j)+\beta}\cap B_j,F_{f(z_j)+\beta}\cap\partial B_j)\to H_{n-1}(F_{f(z_j)+\beta}\cap B_j)\to H_{n-1}(F_{f(z_j)+\beta})$\\
$\to H_{n-1}(F_\alpha)$ if $z_j\not\in X$, the first mapping being the variation.\\
As in the case of Theorem 1.1.4, we may for given $j$ use b) instead of a) if $z_j$ is contained in a stratum of $X$ of positive dimension.\\

{\bf Proof:} By Lemma 1.2.2, we have an epimorphism indeed. The kernel $N$ is generated by the images of 
$$H_n(X_r\cap B_j\cap\{|f-f(z_j)|\le\beta\}, F_{f(z_j)+\beta}\cap B_j)\to H_{n-1}(F_{f(z_j)+\beta})\to H_{n-1}(F_\alpha)$$
see Lemma 1.1.1 and 1.2.1.\\
Finally apply the end of the proof of Theorem 1.1.4 with $f-t_j$ instead of $f$.\\

{\bf Theorem 1.2.4:} Suppose that we have condition $(C_n)$ with $n\ge 1$ or at least $(C''_n)$. Then the kernel in question is generated by the images of the variations\\ 
$H_{n-1}(F_{t_j+\beta},\partial F_{t_j+\beta})\to H_{n-1}(F_{t_j+\beta})\to H_{n-1}(F_\alpha)$.\\ 

{\bf Proof:} The kernel is the image of $H_n(X_r\cap\{|f|\le\alpha\},F_\alpha)\to H_{n-1}(F_\alpha)$. So by Lemma 1.2.1 it is generated by the images of\\
$H_n(X_r\cap\{|f-t_j|\le\beta\},F_{t_j+\beta})\to H_{n-1}(F_{t_j+\beta})\to H_{n-1}(F_\alpha)$,\\ 
i.e. by the kernels of\\ 
$H_{n-1}(F_{t_j+\beta})\to H_{n-1}(X_r\cap\{|f-t_j|\le\beta\})$.\\ 
By Theorem 1.1.5 we finish the proof.\\

The results of this section will be applied in the following section to the local nearby fibre and enable a deeper understanding of local homology groups, they can also be applied to the proof of a global Zariski-van Kampen theorem (section 1.5).\\

{\bf 1.3. First application: hyperplane section of the local nearby fibre}\\

Now we return to the local situation, see section 1.0, with $0\in Z$. We are interested in the fibre $F=F_\alpha$:\\
Let $g:\mathbb{C}^N\to \mathbb{C}$ be a linear function which is sufficiently general. More precisely, it is sufficient to suppose that $g:U\to\mathbb{C}$ is holomorphic and that there is a neighborhood $V$ of $0$ such that $\{g=f=0\}$ is transversal to a Whitney stratification of $(\bar{X},Z)$ in $V\setminus\{0\}$. \\
Choose $0<\alpha\ll\beta\ll\varepsilon\ll 1$. Then $\{\|z\|=\varepsilon\}$ intersects the strata of $\bar{X}\cap \{f=s,g=t\}$ transversally for $|s|\le\alpha,|t|\le\beta\}$. And we have the assumptions of section 1.2 with $F$ instead of $X_r$, $g$ instead of $f$, $\beta$ instead of $\alpha$.\\ 
Now assume $0<\beta'\ll\varepsilon'\ll\alpha\ll\beta\ll\varepsilon$.\\
Let $\{z_1,\ldots,z_s\}$ be the points $z$ of $\bar{F}$ where $\{g=g(z)\}$ is not transversal to $\bar{F}$.\\
Let $B_j:=\{z\in\mathbb{C}^N\,|\,\|z-z_j\|\le \varepsilon'\}$, and let\\
$\{t_1,\ldots,t_k\}:=\{g(z_1),\ldots,g(z_s)\}$ where $k$ is the number of elements of the latter set.\\

First we have:\\

{\bf Lemma 1.3.1:} Let $0<\alpha\ll\beta\ll\varepsilon\ll 1$. Then there is a deformation retraction of $F$ onto $F\cap\{|g|\le\beta\}$ which takes $\partial F$ to $\partial(F\cap\{|g|\le \beta\})$.\\

{\bf Proof:} Let $\phi,\psi:U\to\mathbb{R}$ be defined by $\phi(z)=\|z\|^2, \psi(z)=|g(z)|^2$. By the condition $a$ of Whitney the set of all points $z\in\overline{X}$ with the following condition is open in $\overline{X}$: there is no $\lambda\le 0,c\in \mathbb{C}$ such that\\
$|g(z)|\ge\beta, \|z\|\le\varepsilon, |f(z)|\le\alpha$ and $d\psi(z)|T_zS=\lambda d\phi(z)|T_zS+cdf|T_zS$,\\
where $S$ is the stratum of $\overline{X}$ which contains $z$. By the Curve Selection Lemma we see that this set contains $\{f=0, \|z\|\le\varepsilon, |g(z)|\ge\beta\}$, hence $\{|f|\le\alpha, \|z\|\le\varepsilon, |g(z)|\ge\beta\}$, too. In particular, there is no $z\in S\cap \overline{X}\cap\{f=\alpha\},\lambda\le 0$ such that\\
$|g(z)|\ge\beta, \|z\|\le\varepsilon, d\psi(z)|T_z(S\cap\overline{X}\cap\{f=\alpha\})=\lambda d\phi|T_z(S\cap\overline{X}\cap\{f=\alpha\})$.\\
Therefore we may build a stratified vector field on $\overline{F}$ such that $d\phi(v)>0, d\psi(v)>0$ in the subset $\{\|z\|\le\varepsilon, |g|\ge\beta\}$.\\
Using the flow $\sigma$ of this vector field we can construct the desired homeomorphism: if $z\in X_\varepsilon\cap\{|g|>\beta\}$, there is a $t_z<0$ such that\\
$|g(\sigma(z,t_z))|=\beta$, then $(z,t)\mapsto \sigma(z,tt_z)$, extended by the identity, gives the desired deformation retraction.\\

So, by Lemma 1.3.1 and Theorem 1.2.3, we obtain:\\

{\bf Theorem 1.3.2:} Assume that we have condition $(C_n)$ or at least $(C'_n)$, cf. subsection 1.0. Then the mapping $H_{n-2}(F\cap\{g=\beta\})\to H_{n-2}(F)$ is an epimorphism whose kernel is generated by the images of he mappings:\\
a) $\widetilde{H}_{n-2}(F\cap B_j\cap\{g=g(z_j)+\beta'\})\to H_{n-2}(F\cap\{g=g(z_j)+\beta'\}\to$\\ 
$H_{n-2}(F\cap\{g=\beta\})$ if $z_j\in X$,\\
b) $H_{n-2}(F\cap B_j\cap\{g=g(z_j)+\beta'\},F\cap \partial B_j\cap\{g=g(z_j)+\beta'\})$\\
$\to H_{n-2}(F\cap B_j\cap\{g=g(z_j)+\beta'\})$\\
$\to H_{n-2}(F\cap\{g=\beta\})$ if $z_j\in Z$.\\

Alternatively, by Lemma 1.3.1 and Theorem 1.2.4, we obtain:\\

{\bf Theorem 1.3.3:} Assume that $X$ is smooth and of dimension $n\ge 2$ everywhere or at least that $n=rHd(X)\ge 2$ and $X$ has no strata of dimension $1$. Then the mapping $H_{n-2}(F\cap\{g=\beta\})\to H_{n-2}(F)$ is an epimorphism whose kernel is generated by the images of the mappings\\ 
$H_{n-2}(F\cap\{g=t_j+\beta'\},\partial F\cap \{g=t_j+\beta'\})\to H_{n-2}(F\cap\{g=t_j+\beta'\})\to H_{n-2}(F\cap\{g=\beta\})$,\\
the left one being a variation.\\

{\bf 1.4. Consequences for local homology groups}\\

Let $X_\varepsilon$ be defined as in section 1.0. Furthermore we assume that $X$ is smooth of pure dimension $n$ or at least $rHd(X)=n\ge 2$. We want to discuss consequences for the local homology groups $H_k(X_\varepsilon)$, $k\le n-1$.\\
Let $f:\mathbb{C}^N\to\mathbb{C}$ be a general linear mapping, $H:=f^{-1}(0)$, $F:=X_\varepsilon\cap f^{-1}(\alpha)$, where $0<\alpha\ll\varepsilon$.\\

a) For $k<n-2$ we have
$$H_k(X_\varepsilon)\simeq H_k(X_\varepsilon\cap H)$$
by the (weak) local Lefschetz-Zariski theorem, see [HL2], [HL3].\\
So this case can be reduced to the case $k=n-2$.\\
What is left are the cases $k=n-2, n-1$.\\ 

b) Now look at the case $0\le k=n-2$, hence $n\ge 2$.\\
By the strong local Lefschetz-Zariski theorem (Theorem 1.0.1), we have\\
$H_{n-2}(X_\varepsilon)\simeq H_{n-2}(F)$.\\
Let $g$ be a general linear mapping. 
Then, using Theorem 1.3.2, $H_{n-2}(F)$ can be calculated from $H_{n-2}(F\cap\{g=\beta\})$, where $0<\alpha\ll\beta\ll\varepsilon$.\\

c) Finally look at the case $0\le k=n-1$, hence $n\ge 1$.\\
Here $H_{n-1}(X_\varepsilon)$ can be calculated from $H_{n-1}(F)$, using the local strong Zariski-van Kampen theorem (Cor. 1.0.3).\\

Altogether we can calculate $H_k(X_\varepsilon)$, $0<k\le n-1$, by the $k$-th homology group of a Stein space of dimension $k$.\\

{\bf 1.5. Second application: the global case}\\

The global case of non-singular quasi-projective varieties has been treated by Ch\'eniot [C]. He  uses a blowing-up as usual, we show that one can largely avoid this by reduction to the semi-local case.\\

We drop the non-singularity assumption, as [Ti].\\ 

We proceed as follows. Let $Z\subset \bar{X}$ be subvarieties of $\mathbb{P}_N$, not necessarily irreducible, $X=\bar{X}\setminus Z$ smooth and purely $n$-dimensional, $n>0$. Choose a complex analytic Whitney stratification of $\bar{X}$ which is compatible with $Z$. Let $\{L_t\}$ be a pencil of hyperplanes with base $A$, $A$ being transverse to $\bar{X}$ in the stratified sense. We may assume that $L_t=\{z_1=tz_0\}$, hence $A=\{z_0=z_1=0\}$, and that $L_0$ and $L_\infty$ are general. Let $p_1,\ldots,p_s$ be the points in which the corresponding hyperplane $L_{\tau_j}$ through $p_j$ is not transversal to $\bar{X}$. Choose $\alpha>0$ such that\\
$|\tau_1|<\alpha,\ldots,|\tau_s|<\alpha$, and put $D=\{|t|\le\alpha\}$. Let $\{t_1,\ldots,t_k\}:=\{\tau_1,\ldots,\tau_s\}$ if the latter set has $k$ elements. Let $0<\beta\ll\varepsilon\ll\alpha$, $F_t:=X\cap L_t$. On $\mathbb{P}_N\setminus L_\infty$ we take inhomogeneous coordinates $\xi_1:=z_1/z_0,\ldots, \xi_N:=z_N/z_0$. For $j\in\{1,\ldots,s\}$ let $B_j=\{\|\xi-p_j\|\le\varepsilon\}$ be the closed ball around $p_j$ of radius $\varepsilon$.\\
Let $f:X\setminus A\to\mathbb{P}_1$ be defined by $f=\frac{z_1}{z_0}$. Then $F_t$ is the closure of $f^{-1}(t)$ in $X$.\\
For $M\subset\mathbb{P}_1$ let $X_M$ be the closure of $f^{-1}(M)\cap X$ in $X$.\\
Let $r\gg 0$ and ${\rm B}_r:=\{\|\xi\|\le r\}$.\\
First:\\

{\bf Lemma 1.5.1:} The mapping $f:\bar{X}\cap\partial {\rm B}_r\to \mathbb{C}$ has no critical points which are mapped to $D$.\\

{\bf Proof:} Let $q_{\nu}$ be a sequence of points in $\mathbb{P}_N\setminus L_\infty$ which converges to a point $q$ in $L_\infty$. Let $r_\nu:=\|q_\nu\|$ and $T_\nu$ the tangent space to $\partial{\rm B}_r$ at $q_\nu$. Assume that $T_\nu$ converges to a subspace $T$ of the tangent space of $\mathbb{P}_N$ at $q$. Then it is easy to see that the tangent space of $L_\infty$ at $q$ is contained in $T$.\\
By hypothesis, all $L_t$, $t\in D$, intersect $L_\infty$ transversally, so $\partial {\rm B}_r$, too. The same holds if we intersect by $\bar{X}$. This implies our statement.\\

So we can apply section 1.3. Suppose that $X$ is smooth of pure dimension $n\ge 1$ or at least that $n=rHd(X)$.\\

{\bf Lemma 1.5.2:} a) $H_k(X,X_D)=0, k\le n$.\\
b) For any $k$, we have $H_k(X_D\cap {\rm B}_r,F_\alpha\cap {\rm B_r})\simeq H_k(X_D,F_\alpha)$.\\

{\bf Proof:} b) By excision, $H_k(X_D\cap{\rm B}_r, F_\alpha\cap{\rm B}_r)\simeq H_k((X_D\cap{\rm B}_r)\cup F_\alpha, F_\alpha)$.\\
Now $(X_D\cap{\rm B}_r)\cup F_\alpha$ is a weak deformation retract of $X_D$:\\
It is sufficient to show that $(X_D\cap\partial{\rm B}_r)\cup(F_\alpha\setminus\stackrel{o}{\rm B}_r)$ is a weak deformation retract of $X_D\setminus\stackrel{o}{\rm B}_r$.\\
First let us show that $X\cap A$ is a deformation retract of both $F_\alpha\setminus\stackrel{o}{\rm B}_r$ and $X_D\setminus\stackrel{o}{\rm B}_r$: \hfill{(*)}\\
Let $\widetilde{X}\subset \mathbb{C}^{N+1}\setminus\{0\}$ be the cone over $X$ and $\pi:\widetilde{X}\to X$ the projection,\\
$\widetilde{X}_\varepsilon:=\widetilde{X}\cap\{\|z\|\le\varepsilon\}$. Then $X$ is the orbit space of $\partial\widetilde{X}_\varepsilon$ with respect to the usual $S_1$ action. And $\{z_0=z_1=0\}$ intersects $\partial\widetilde{X}_\varepsilon$ transversally, so $\{z_1=t _1,z_2=t_2\}$, too, if $t_1,t_2$ are sufficiently small.\\
Now $\{(t_1,t_2)\in\mathbb{C}^2\,|\,t_1=\alpha t_0,|t_0|\le\beta\}$ is contractible, where $0<\beta\ll 1$. Then $\partial\widetilde{X}_\varepsilon\cap\{z_0=z_1=0\}$ is a deformation retract of\\
$\partial\widetilde{X}_\varepsilon\cap\{z_1=\alpha z_0, |z_0|\le\beta\}$. As  for the images under $\pi$, the image of\\
$\partial\widetilde{X}_\varepsilon\cap\{z_0=z_1=0\}$ is $X\cap A$,\\
the image of $\partial\widetilde{X}_\varepsilon\cap\{z_1=\alpha z_0, 0<|z_0|\le\beta\}$ is $F_\alpha\setminus (L_\infty\cup\stackrel{o}{\rm B}_r)$, with\\ $r=\sqrt{\varepsilon^2/\beta^2-1}$:\\
In fact, $\partial\widetilde{X}_\epsilon$ is given in $\widetilde{X}$ by the equation $|z_0|^2+\ldots+|z_N|^2=\varepsilon^2$, the condition $0<|z_0|\le\beta$ leads to 
$$z_0\neq 0\,,1+\frac{|z_1|^2}{|z_0|^2}+\ldots+\frac{|z_N|^2}{|z_0|^2}=\frac{\varepsilon^2}{|z_0|^2}\ge\frac{\varepsilon^2}{\beta^2}$$
hence 
$$z_0\neq 0\,, \frac{|z_1|^2}{|z_0|^2}+\ldots+\frac{|z_N|^2}{|z_0|^2}\ge r^2$$ 
with $r^2=\frac{\varepsilon^2}{\beta^2}-1$.\\
In $X\setminus L_\infty$, this corresponds to $\xi\not\in \stackrel{o}{\rm B}_r$.\\
So we have a commutative diagram
$$\begin{array}{ccc}
\partial\widetilde{X}_\varepsilon\cap\{z_0=z_1=0\}&\to&\partial\widetilde{X}_\varepsilon\cap\{z_1=\alpha z_0,|z_0|\le\beta\}\\
\downarrow&&\downarrow\\
X\cap A&\to&F_\alpha\setminus \stackrel{o}{\rm B}_r
\end{array}$$
We have seen that the upper horizontal gives a (weak) deformation retract. Since the verticals correspond to $S_1$-fibrations we obtain that the lower horizontal gives a (weak) deformation retract, too. This can be seen by the property of relative homotopy groups for fibrations, with a little extra care in case of $\pi_0$.\\
This proves the first half of the statement (*).\\
In an analogous way we prove the second half, replacing the condition\\
$z_1=\alpha z_0$ by $|z_1|\le\alpha |z_0|$ and using the commutative diagram
$$\begin{array}{ccc}
\partial\widetilde{X}_\varepsilon\cap\{z_0=z_1=0\}&\to&\partial\widetilde{X}_\varepsilon\cap\{|z_1|\le\alpha |z_0|,|z_0|\le\beta\}\\
\downarrow&&\downarrow\\
X\cap A&\to&X_D\setminus \stackrel{o}{\rm B}_r
\end{array}$$
Altogether this shows $(*)$.\\
Furthermore, note that $F_\alpha\cap\partial{\rm B}_r$ is a deformation retract of $X_D\cap\partial{\rm B}_r$, because of Lemma 1.5.1. By excision, we get that $F_\alpha\setminus \stackrel{o}{\rm B}_r$ is a deformation retract of\\
$(X_D\cap\partial{\rm B}_r)\cup(F_\alpha\setminus\stackrel{o}{\rm B}_r)$. Now (*) implies that $F_\alpha\setminus\stackrel{o}{\rm B}_r$ is a deformation retract of  $X_D\setminus\stackrel{o}{\rm B}_r$, so we conclude that $(X_D\cap\partial{\rm B}_r)\cup(F_\alpha\setminus\stackrel{o}{\rm B}_r)$ is a deformation retract of $X_D\setminus\stackrel{o}{\rm B}_r$. Altogether this proves our statement that $(X_D\cap\partial{\rm B}_r)\cup(F_\alpha\setminus\stackrel{o}{\rm B}_r)$ is a weak deformation retract of $X_D\setminus\stackrel{o}{\rm B}_r$.
\vskip.1in
a) By excision it is sufficient to show that\\ $H_k(X_{\mathbb{P}_1\setminus\stackrel{o}{D}},X_{\partial D})=0, k\le n$.\\ 
Using $\mathbb{P}_1\to\mathbb{P}_1: t\mapsto 1/t$ we reduce to a statement of the form\\ $H_k(X_{D'},X_{\partial D'})=0, k\le n$, with $\alpha'>0$ sufficiently small, $D':=\{|t|\le\alpha'\}$.\\
Nów $f:X^r_{D'}\cap{\rm B}_r\to D'$ defines a trivial fibre bundle, so we get that\\
$(X^r_{D'}\cap {\rm B}_r,\partial X^r_{D'}\cap {\rm B}_r) $ is homeomorphic to $(F_{\alpha'}\cap{\rm B}_r,F_{\alpha'}\cap\partial {\rm B}_r)\times({\mathbb D}^2,\partial{\mathbb D}^2)$.\\
And $H_k(F_{\alpha'}\cap{\rm B}_r,F_{\alpha'}\cap\partial {\rm B}_r)=0, k\le n-2$,\\
because of stratified Morse theory, applied to $z\mapsto -\|z\|^2$,\\
so $H_k(X_{D'}\cap{\rm B}_r,\partial X_{D'}\cap{\rm B}_r)=0, k\le n$.\\
By b) we have that\\
$H_k(X_{D'}\cap {\rm B}_r,F_{\alpha'}\cap {\rm B_r})\simeq H_k(X_{D'},F_{\alpha'})$.\\
In the same way we can show that\\
$H_k(X_{\partial D'}\cap {\rm B}_r,F_{\alpha'}\cap {\rm B_r})\simeq H_k(X_{\partial D'},F_{\alpha'})$:\\
 instead of the inequality $|z_1|\le\alpha' |z_0|\le\beta$ use the inequality $|z_1|=\alpha' |z_0|\le\beta$.\\
So we get\\
$0=H_k(X_{D'}\cap{\rm B}_r,X_{\partial D'}\cap {\rm B}_r)\simeq H_k(X_{D'},X_{\partial D'})$
for $k\le n$.\\

{\bf Theorem 1.5.3:} a) $H_k(X,F_\alpha)=0$ for $k<n$, so $H_{n-1}(F_\alpha)\to H_{n-1}(X)$ is surjective,\\ 
b) the kernel of the epimorphism $H_{n-1}(F_\alpha)\to H_{n-1}(X)$ is $N_1+\ldots+N_s$, where $N_j$ is the image of\\
$\widetilde{H}_{n-1}(F_{\tau_j+\beta}\cap B_j)\to H_{n-1}(F_{\tau_j+\beta})\to H_{n-1}(F_\alpha)$ if $p_j\in X$,\\
$H_{n-1}(F_{\tau_j+\beta}\cap B_j,F_{\tau_j+\beta}\cap\partial B_j)\to H_{n-1}(F_{\tau_j+\beta}\cap B_j)\to H_{n-1}(F_\alpha)$ if $p_j\in Z$, the first mapping given by the local variation.\\ 

{\bf Proof:} a) For $k<n$: $H_k(X_D\cap {\rm B}_r,F_\alpha\cap{\rm B}_r)=0$, by Lemma 1.2.2, hence $H_k(X_D,F_\alpha)=0$, by Lemma 1.5.2.\\
b) Look at the commutative diagram
$$\begin{array}{ccc}
H_n(X_D\cap {\rm B}_r,F_\alpha\cap {\rm B}_r)&\to&H_{n-1}(F_\alpha\cap {\rm B}_r)\\
\downarrow&&\downarrow\\
H_n(X,F_\alpha)&\to&H_{n-1}(F_\alpha)
\end{array}$$
The left vertical is surjective: $H_n(X_D\cap {\rm B}_r,F_\alpha\cap {\rm B}_r)\to H_n(X_D,F_\alpha)$ is bijective by Lemma 1.5.2b), and $H_n(X_D,F_\alpha)\to H_n(X,F_\alpha)$ is surjective by Lemma 1.5.2a).\\
So the kernel in question is the image of the diagonal mapping $H_n(X_D\cap {\rm B}_r,F_\alpha\cap {\rm B}_r)\to H_{n-1}(F_\alpha)$.\\
Now the image of the upper horizontal has been described in Theorem 1.2.3, which implies our statement.\\

Under more restrictive conditions we have:\\

{\bf Theorem 1.5.4} (see [C]): Let us assume that $X$ is smooth of dimension $n$ (or at least that $rHd(X)\ge n$ and $X$ has no strata of dimension $0$). Then the kernel of the epimorphism $H_{n-1}(F))\to H_{n-1}(X)$ is the image of the mappings\\
$H_{n-1}(F,F\cap A)\to H_{n-1}(F_{t_j+\beta},F\cap A)\to H_{n-1}(F_{t_j+\beta})\to H_{n-1}(F)$\\ 
where the second mapping is defined by $[c]\mapsto [(h_j)_*(c)-c]$, $h_j$ being a geometric monodromy.\\

{\bf Proof:} The assumption allows to assume that the points of non-transversality are contained in $Z$. Furthermore argue as in the case of Theorem 1.2.4.\\

{\bf 1.6. Additional comments}\\

Here we concentrate upon the computation of the first homology group, with $n\ge 2$. In special cases there are alternatives:\\

a) Take the situation of section 1.0 with $\bar{X}=U, Z$ hypersurface in $U$, $n=N$:\\
Then $H_1(X_\varepsilon)=H_1(B_\varepsilon\setminus Z)\simeq H_1(S_\varepsilon\setminus Z)\simeq H^{2n-2}(S_\varepsilon,S_\varepsilon\cap Z)\simeq$\\
$H^{2n-3}(S_\varepsilon\cap Z)\simeq H^{2n-3}(S_\varepsilon\cap Z,S_\varepsilon\cap Sing\,Z)\simeq H_0(S_\varepsilon\cap(Z\setminus Sing\,Z))\simeq \mathbb{Z}^k$,\\ 
where $k$ is the number of irreducible components of $Z$.\\
In the special case $n=2$ this is well-known from knot theory: then $S_\varepsilon\cap Z$ is the disjoint union of $k$ embedded circles in $S_3$, so we have a knot if $k=1$, otherwise a link.\\

b) Similarly in the global case where $\bar{X}=\mathbb{P}_n(\mathbb{C})$ and $Z$ a hypersurface in $\mathbb{P}_n$.\\
Then $H_1(\mathbb{P}_n\setminus Z)\simeq H^{2n-1}(\mathbb{P}_n,Z)$.\\
Suppose that $Z$ has $k$ irreducible components of degrees $d_1,\ldots,d_k$. Let $d$ be the greatest common divisor of $d_1,\ldots,d_k$. Then the mapping\\
$H^{2n-2}(\mathbb{P}_n)\to H^{2n-2}(Z)$ corresponds to the homomorphism $\mathbb{Z}\to \mathbb{Z}^k$\\
with $1\mapsto (d_1,\ldots, d_k)$, so we have a short exact sequence
$$0\to H^{2n-2}(\mathbb{P}_n)\to H^{2n-2}(Z)\to H^{2n-1}(\mathbb{P}_n,Z)\to 0$$
which implies\\
$H_1(\mathbb{P}_n\setminus Z)\simeq H^{2n-1}(\mathbb{P}_n,Z)\simeq (\mathbb{Z}/d\mathbb{Z})\times\mathbb{Z}^{k-1}$,\\
where $d$ is the greatest  common divisor of $d_1,\ldots,d_k$.\\

\centerline{\bf 2. Homotopy groups}
\vskip.1in

{\bf 2.0. Preparation}\\

We try to transpose our methods from homology to homotopy groups.\\

The main difficulty comes from excision. Suppose that $A_1,A_2$ as well as their union and intersection are $CW$-complexes.\\
By excision, we have that $H_k(A_1,A_1\cap A_2)\simeq H_k(A_1\cup A_2,A_2)$ for all $k$.\\
This is no longer true if we pass to homotopy groups, we merely have the\\

{\bf Blakers-Massey theorem:}\\
{\it Suppose that $(A_1,A_1\cap A_2)$ is $m$-connected and that $(A_2,A_1\cap A_2)$ is $n$-connected, $m\ge 1,n\ge 0$, $\cdot\in A_1\cap A_2$. Then the mapping\\
$\pi_k(A_1,A_1\cap A_2,\cdot)\to \pi_k(A_1\cup A_2,A_2,\cdot)$ is bijective for $1\le k<m+n$ and surjective for $1\le k=m+n$.} \\
See [tD].\\

Now let $A$ be a CW complex, $B_1,B_2,C$ subcomplexes, $A=B_1\cup B_2, C=B_1\cap B_2$.
Then the Mayer-Vietoris sequence implies that $H_k(B_1,C)\oplus H_k(B_2,C)\simeq H_k(A,C)$. What happens if we pass to homotopy groups?\\

Assume that $(B_i,C)$ is $(n-1)$-connected, $i=1,2$, where $n\ge 2$, $\cdot\in C$.\\

{\bf Lemma 2.0.1:} a) $(A,C)$ is $(n-1)$-connected, too.\\
b) $\pi_n(B_1,C,\cdot)\oplus \pi_n(B_2,C,\cdot)\to\pi_n(A,C,\cdot)$ is surjective if $n\ge 3$,\\
$\pi_2(B_1,C,\cdot)\star \pi_2(B_2,C,\cdot)\to\pi_2(A,C,\cdot)$ is surjective if $n=2$.\\
Here $\star$ denotes the free product.\\

{\bf Proof:} By Blakers-Massey, $\pi_k(B_1,C,\cdot)\to \pi_k(A,B_2,\cdot)$ is bijective for \\
$0<k<n$ and surjective for $k=n$ because $n\le 2n-2$.\\
We have an exact sequence $\pi_k(B_2,C,\cdot)\to \pi_k(A,C,\cdot)\to \pi_k(A,B_2,\cdot)$, $k>0$, by the exact homotopy sequence for a triple, see [Sp] Theorem 7.2.15.
\vskip.1in
a) It is clear that each point of $A$ can be joined to a point of $C$, so $(A,C)$ is $0$-connected.\\
If $0<k<n$, $\pi_k(B_2,C,\cdot)=0$ and $\pi_k(A,B_2,\cdot)=\pi_k(B_1,C,\cdot)=0$, see beginning of the proof, so $\pi_k(A,C,\cdot)=0$.
\vskip.1in
b) Now let $x\in \pi_n(A,C,\cdot)$. Look at the following commutative diagram:
$$\begin{array}{ccccc}
&&\pi_n(B_1,C,\cdot)&&\\
&&\downarrow&\searrow&\\
\pi_n(B_2,C,\cdot)&\to&\pi_n(A,C,\cdot)&\to&\pi_n(A,B_2,\cdot)
\end{array}$$
The image of $x$ in $\pi_n(A,B_2,\cdot)$ comes from an element $y_1\in\pi_n(B_1,C,\cdot)$ (because $\pi_n(B_1,C,\cdot)\to \pi_n(A,B_2,\cdot)$ is surjective, see beginning of the proof) which is mapped to an element $x_1\in\pi_n(A,C,\cdot)$, so $x$ and $x_1$ are mapped to the same element of $\pi_n(A,B_2,\cdot)$. Hence $x_1^{-1}x$ is mapped to $1\in\pi_n(A,B_2,\cdot)$.\\ 
So $x_1^{-1}x$ comes from an element $y_2\in \pi_n(B_2,C,\cdot)$. If $n\ge 3$, this implies that $\pi_n(B_1,C,\cdot)\oplus \pi_n(B_2,C,\cdot)\to\pi_n(A,C,\cdot)$ is surjective because $(y_1,y_2)$ is mapped to $x$. Similarly in the case $n=2$ where we dispose only of a homomorphism 
$\pi_2(B_1,C,\cdot)\star \pi_2(B_2,C,\cdot)\to\pi_2(A,C,\cdot)$.\\

Now let us return to the data of section 1.0.\\

First recall the condition\\

$(C_n)$: $X$ is smooth, $\dim(X,x)\ge n$ for all $x\in X$.\\

We modify the conditions $C'_n$ and $C''_n$: Note that {\bf $(C_n)$} implies the condition\\

$({\cal C}'_n)$: $rhd(X)\ge n$, where $rhd$ means the rectified {\it homotopical} depth.\\ 
This means that for every stratum $S$ with $dim\,S=k$ and every $x\in S$ there is a fundamental system of neighborhoods $V$ of $x$ in $X$ such that $(V,V\setminus S)$ is $n-k-1$-connected.\\

Note that ${\cal C}'_n$ holds if $X$ is a local complete intersection with $dim(X,x)\ge n$ for all $x\in X$, and ${\cal C}'_n$ implies  $C'_n$.\\

If $n\ge 1$, $(C_n)$ implies the condition\\

$({\cal C}''_n)$: As in the case of $({\cal C}'_n)$, with $(V,V\setminus S)$ being even $n$-connected if $dim\,S=0$.\\

And $({\cal C}'_n)$ implies $({\cal C}''_n)$ if all strata are of dimension $>0$, of course. Furthermore, $({\cal C}'_n)$ implies $(C'_n)$, $({\cal C}''_n)$ implies $(C''_n)$.\\

{\bf Lemma 2.0.2:} a) Suppose that $rhd(X)\ge n$. Then for every $x\in S$, $S$ being a stratum of dimension $k$, there is a fundamental system of neighborhoods $V$ of $x$ in $X$ such that $(V,V\setminus\{x\})$ is $(n+k-1)$-connected.\\
b) $({\cal C}''_n)$ implies that for every $x\in X$ there is a fundamental system of neighborhoods $V$ of $x$ in $X$ such that $(V,V\setminus\{x\})$ is $n$-connected.\\

{\bf Proof:} a) $(V,V\setminus\{x\})$ may be replaced by $(V,\partial V)$, $V$ being suitably chosen, which is homeomorphic to $(N,\partial N)\times (I^{2k},\partial I^{2k})$, where $N$ is a normal slice and $I:=[0,1]$. Now $(N,\partial N)$ has the same homotopy type as $(V,V\setminus S)$, so it is $(n-k-1)$-connected. By [H2] Lemma 1.1c) we have that\\
$(N,\partial N)\times (I^{2k},\partial I^{2k})$ is $(n+k-1)$-connected.\\
b) follows from a).\\
 
Now we have a version of the strong local Lefschetz-Zariski theorem for homotopy groups:\\

{\bf Theorem 2.0.3 (Strong local Lefschetz-Zariski theorem} ([HL2], [HL3]):\\
The pair $(X_\varepsilon,F)$ is $(n-1)$-connected, hence $\pi_k(X_\varepsilon,F,\cdot)=0, 1\le k<n$, 
and\\
$\pi_k(F,\cdot) \to \pi_k(X_\varepsilon,\cdot)$ is bijective for $k<n-1$, surjective for $k=n-1$.\\
Here, $\cdot\in F$ is a base point.\\

Similarly, we derived in [H2] a strong local Zariski-van Kampen theorem for homotopy groups:
Let $W_k: \pi_k(F,\partial F,\cdot)\to \pi_{k+1}(X_\varepsilon,F,\cdot)$ be the composition of the following mappings, with $C:=\partial X_\varepsilon\cap\{|f|=\alpha,Im\,f\ge 0\}$:
$$\begin{array}{c}
\pi_k(F,\partial F,\cdot)\\
\downarrow\\
\pi_k((\{f=\pm\alpha\}\cap X_\varepsilon)\cup C,(\{f=-\alpha\}\cap X_\varepsilon)\cup C,\cdot)\\
\simeq\uparrow\\
\pi_{k+1}(\{|f|=\alpha,Im\,f\ge 0\}\cap X_\varepsilon,(\{f=\pm\alpha\}\cap X_\varepsilon)\cup C,\cdot)\\
\downarrow\\
\pi_{k+1}(\{|f|\le\alpha\}\cap X_\varepsilon),(\{f=\pm\alpha\}\cap X_\varepsilon)\cup C,\cdot)\\
\downarrow\\
\pi_{k+1}(\partial(X_\varepsilon\cap\{|f|\le\alpha\}),(\{|f|=\alpha,Im\,f\ge 0\}\cap X_\varepsilon)\cup\partial(X_\varepsilon\cap\{|f|\le\alpha\}),\cdot)\\
\simeq\uparrow\\
\pi_{k+1}(\partial(X_\varepsilon\cap\{|f|\le\alpha\},F,\cdot)\\
\simeq\downarrow\\
\pi_{k+1}(X_\varepsilon,F,\cdot)
\end{array}$$ 
Then we have:\\

{\bf Theorem 2.0.4:} Assume $n\ge 2$, $(F,\partial F)$ $(n-2)$-connected.\\
a) $(X_\varepsilon,F)$ is $(n-1)$-connected, and $W_{n-1}:\pi_{n-1}(F,\partial F,\cdot)\to \pi_n(X_\varepsilon,F,\cdot)$ is bijective for $n\ge 4$, surjective for $n=3$, and generates the target group for $n=2$.\\
b) $\partial\circ W_{n-1}=\kappa\circ Var_{n-1}(\frac{f}{\alpha})$, where $Var_{n-1}(\frac{f}{\alpha})(\beta)=\beta\bot (h\circ\beta)$ and $\kappa$ is the inversion.\\
c) $\pi_{n-1}(F,\cdot)/N\simeq \pi_{n-1}(X_\varepsilon,\cdot)$ where $N$ is the image of $Var_{n-1}(\frac{f}{\alpha})$ if $n\ge 3$ resp. the subgroup generated by it if $n=2$.\\

{\bf Proof:} with [H2] Theorem 6.6. And c) is already derived in [H2] Theorem 0.2 if $X$ is smooth.\\

As for calculation, note that we may compute $Im\,Var_{n-1}(f)$ as soon as we have generators of $\pi_{n-1}(F,\partial F,\cdot)$ if $n\ge 3$ because $Var_{n-1}(f)$ is a homomorphism in this case. The case $n=2$ is exceptional. It is no restriction to suppose that $X_\varepsilon$ and hence $F$ is connected. Then:\\

{\bf Theorem 2.0.5:}
Suppose that $n=2$, $F$ connected. Let $a_1,\ldots,a_r\in\pi_1(F,\partial F,\cdot)$ be the images of generators of $\pi_1(F,\cdot)$. Let $z_0,\ldots,z_s$ be representatives of the connected components of $\partial F$, where $z_0=\cdot$, and let $b_j\in \pi_1(F,\partial F,\cdot)$ be represented by a path from $z_j$ to $z_0=\cdot$. Then the subgroup $U$ of $\pi_1(F,\cdot)$ generated by $Im\,Var_1(f)$ coincides with the normal subgroup $N$ of $\pi_1(F,\cdot)$ generated by the images of $a_1,\ldots,a_r,b_1,\ldots,b_s$ under $Var_1(f)$.\\

{\bf Proof:}  $U$ is a normal subgroup of $\pi_1(F,\cdot)$  because it coincides with\\ $ker(\pi_1(F,\cdot)\to \pi_1(X_\varepsilon,\cdot))$. Of course, $N\subset U$. So it remains to show that $U\subset N$.\\
Note that $\pi_1(F,\{z_0,\ldots,z_s\},\cdot)\to \pi_1(F,\partial F,\cdot)$ is surjective, whereas\\ 
$\pi_1(F,\cdot)\to \pi_1(F,\{z_0,\ldots,z_s\},\cdot)$ is injective, so we consider $\pi_1(F,\cdot)$ as a subset of $\pi_1(F,\{z_0,\ldots,z_s\},\cdot)$. So we may define $Var_1(f)$ on $\pi_1(F,\{z_0,\ldots,z_s\},\cdot)$ instead and accomplish the calculation inside $\pi_1(F,\{z_0,\ldots,z_s\},\cdot)$ then.\\ 
Now: $Var_1(f)(ba^{-1})\in N$ if\\
$b\in \pi_1(F,\{z_0,\ldots,z_s\},\cdot),a\in \pi_1(F,\cdot), Var_1(f)(a), Var_1(f)(b)\in N$.\\
This is because $Var_1(f)(ba^{-1})=ab^{-1}h(b)h(a)^{-1}=a(b^{-1}h(b)h(a)^{-1}a)a^{-1}
$\\
$=a\,Var_1(f)(b)(Var_1(f)(a))^{-1})a^{-1}\in N$.\\
Therefore $Var_1(f)(a)\in N$ if $a\in \pi_1(F,\cdot)$, so $Var_1(f)(b_ja)\in N$, hence $U\subset N$: note that the elements of  $\pi_1(F,\{z_0,\ldots,z_s\},\cdot)$ are of the form $b_ja, a\in \pi_1(F,\cdot)$.\\

In the following we pass to the semi-local case.\\

{\bf 2.1. Towards the semi-local case}\\

Now we would like to pass from the local to a semi-local situation, first using the setup of section 1.1. \\
We assume condition $C_n$ or at least ${\cal C}'_n$, where $n\ge 2$. Let $\cdot\in F_\beta$.\\

{\bf Lemma 2.1.1:} a) We have an epimorphism 
$$\star_{j=1}^k \pi_n(F_\beta\cup(X\cap B_j\cap\{|f|\le\beta\}),F_\beta,\cdot) \to
\pi_n(X_r\cap\{|f|\le\beta\},F_\beta,\cdot)$$
If $n\ge 3$ we may replace $\star$ by $\oplus$.\\
b) $\pi_{n-1}(F_\beta,\cdot)\to \pi_{n-1}(X_r\cap\{|f|\le\beta\},\cdot)$ is an epimorphism, the kernel being generated by the images of $\pi_n(F_\beta\cup(X\cap B_j\cap\{|f|\le\beta\}),F_\beta,\cdot) \to\pi_{n-1}(F_\beta,\cdot)$.\\

{\bf Proof:} a) We may replace the right hand side by $\pi_n(F_\beta\cup(\bigcup_{j=1}^k X\cap B_j\cap\{|f|\le\beta\}),F_\beta,\cdot)$, see proof of Lemma 1.1.1. Then apply Lemma 2.0.1b).\\
b) The pairs $(X\cap B_j\cap\{|f|\le\beta\},F_\beta\cap B_j)$ are $(n-1)$-connected, because of local strong Lefschetz-Zariski, hence $(F_\beta\cup(\bigcup_{j=1}^k B_j\cap X\cap\{|f|\le\beta\}),F_\beta)$, too, the first space being a deformation retract of $X_r\cap\{|f|\le\beta\}$, see proof of Lemma 1.1.1. So we have an epimorphism.\\
The kernel is the image of $\pi_n(X_r\cap\{|f|\le\beta\},F_\beta,\cdot)$.\\
The rest follows from a).\\

Now modify the hypothesis using ${\cal C}''_n$ instead of ${\cal C}'_n$.\\

Then we have an analogue of Theorem 1.1.4 and 1.1.5. Let $r>0$ be general, $0<\beta\ll\varepsilon\ll r$:\\

{\bf Theorem 2.1.2:} Let $N$ be the kernel of the epimorphism $\pi_{n-1}(F_\beta,\cdot)\to \pi_{n-1}(X_r\cap\{|f|\le\beta\},\cdot)$. Then\\
a) $N$ is the subgroup generated by $N_1,\ldots,N_k$, where $N_j$ is the subgroup generated by
the image of the variation mapping\\
$\pi_{n-1}(F_\beta,F_\beta\setminus\stackrel{o}{B}_j,\cdot)\to \pi_{n-1}(F_\beta,\cdot)$ defined by $[\gamma]\mapsto [\gamma\bot(h_j\circ\gamma)]$,\\
where $h_j:F_\beta\to F_\beta$ is a geometric monodromy of $f|X\cap B_j\cap\{|f|=\beta\}$, extended by the identity,\\
b) $N$ is the subgroup generated by the image of the variation mapping $\pi_{n-1}(F_\beta,\partial F_\beta,\cdot)\to \pi_{n-1}(F_\beta,\cdot)$ of $f:X_r\cap\{|f|=\beta\}\to \beta S_1$, defined by $[\gamma]\mapsto [\gamma\bot(h\circ\gamma)]$, where $h$ is a geometric monodromy of\\
 $f:X_r\cap\{|f|=\beta\}\to \beta S_1$.\\

{\bf Proof:} We have an epimorphism, see Lemma 2.1.1b), and may assume that $z_1,\ldots,z_k\in Z$ because we may replace $X$ by $X\setminus\{z_1,\ldots,z_k\}$: for every $z_j$ there is a fundamental system of neighbourhoods $V$ in $X$ such that $(V,V\setminus \{z_j\}$ is $n$-connected, see Remark after Lemma 2.0.1.\\
a) By Lemma 2.1.1a), $N$ is generated by $\widetilde{N}_1,\ldots,\widetilde{N}_k$, where $\widetilde{N}_j$ is generated by the image of $\pi_n(F_\beta\cup(X\cap B_j\cap\{|f|\le\beta\}),F_\beta,\cdot) \to\pi_{n-1}(F_\beta,\cdot)$. Fix $j$ and a trivialization $\phi_j:(\partial B_j\cap F_\beta)\times \beta\mathbb{D}_2\to X\cap\partial B_j\cap\{|f|\le\beta\}$ of $f|X\cap\partial B_j\cap\{|f|\le\beta\}$ such that $\phi(z,\beta)=z$. Use $\phi_j$ in order to attach $(F_\beta\setminus \stackrel{o}{B_j})\times\beta\mathbb{D}_2$ to $X\cap B_j\cap\{|f|\le\beta\}$, we obtain a space called $X_j$. Then $(X_j\cap\{|f|=\beta\})\cup(X_j\setminus \stackrel{o}{B}_j)$ defines an open book. The rest follows from [H2] Theorem 6.6.\\
b) Note that the space $X_r\cap(\{|f|=\beta\}\cup\{|f|\le\beta\}\setminus\bigcup_{j=1}^k\stackrel{o}{B}_j)$ is a deformation retract of $X_r\cap\{|f|\le\beta\}$ and defines an open book. So we can use the corresponding variation which uses some monodromy which is the identity on $F_\beta\setminus\bigcup_{j=1}^k\stackrel{o}{B}_j$.\\
Now the mapping $\pi_{n-1}(F_\beta,\partial F_\beta)\to \pi_{n-1}(F_\beta,F_\beta\setminus \bigcup_{j=1}^k\stackrel{o}{B}_j)$ is surjective:\\
First, the pair $(F_\beta\setminus\bigcup\stackrel{o}{B}_j,\partial F_\beta)$ is $(n-2)$-connected: This follows by stratified Morse theory, using the function $z\mapsto -\|z\|^2$, see proof of Theorem 1.1.5.\\
Then, for $n>2$, conclude using the long exact homotopy sequence of the triple 
$(F_\beta,F_\beta\setminus\bigcup\stackrel{o}{B}_j,\partial F_\beta)$.\\
In the case $n=2$ the desired surjectivity follows from the fact that\\ 
$(F_\beta\setminus\bigcup\stackrel{o}{B}_j,\partial F_\beta)$ is $0$-connected.\\
Note that, as a result, this means that we can take the variation associated with the open book $X_r\cap\{|f|=\beta\}\cup\partial X_r\cap\{|f|\le\beta\}$ - which was not clear from the beginning because this open book need not be a deformation retract of $X_r\cap\{|f|\le \beta\}$.\\

{\bf 2.2. The semi-local case}\\

Now turn to the setup of section 1.2. Let $\Gamma_j$ be the trace of a suitable path $\gamma_j$ from $\alpha$ to $t_j+\beta$, see Theorem 1.2.1. Choose a base point $\cdot$ in $F_\alpha$. We define $\rho_j$ as a lifting of $\gamma_j$ with $\rho_j(0)=\cdot$  and put $\cdot_j:=\rho_j(1)$.\\
We assume that condition $C_n$ or at least ${\cal C}'_n$ holds, with $n\ge 2$.\\

{\bf Lemma 2.2.1:} We have an epimorphism
$$\star_{j=1}^k\pi_n(X_r\cap \{|f-t_j|\le\beta\},F_{t_j+\beta},\cdot_j)\to \pi_n(X_r\cap\{|f|\le\alpha\},F_\alpha,\cdot)$$
If $n\ge 3$ we may replace $\star$ by $\oplus$.\\

{\bf Proof:} Note that $(X_r\cap \{|f-t_j|\le\beta\},F_{t_j+\beta})$ can be replaced above by $(X_r\cap f^{-1}(\Delta_j),X_r\cap f^{-1}(\Gamma_j))$. Then, $\cdot_j$ can be replaced by $\cdot$. So we have passed to $\star_{j=1}^k\pi_n(X_r\cap f^{-1}(\Delta_j),X_r\cap f^{-1}(\Gamma_j),\cdot)$ which is mapped surjectively to $\pi_n(X_r\cap f^{-1}(\Delta),X_r\cap f^{-1}(\Gamma),\cdot))$, by Lemma 2.0.1. Finally, as remarked in the proof of Lemma 1.2.1, $X_r\cap f^{-1}(\Delta)$ is a deformation retract of $X_r\cap\{|f|\le\alpha\}$, and $F_\alpha$ is a deformation retract of $X_r\cap f^{-1}(\Gamma)$.\\

Furthermore we assume that $(C_n)$ or at least $({\cal C}''_n$) holds, with $n\ge 2$.\\

{\bf Theorem 2.2.2 (Semi-local Zariski-van Kampen theorem for homotopy groups):} The mapping $\pi_{n-1}(F_\alpha,\cdot)\to \pi_{n-1}(X_r\cap\{|f|\le\alpha\},\cdot)$ is an epimorphism whose kernel is generated by the images of the variation mappings\\
$\pi_{n-1}(F_{f(z_j)+\beta},F_{f(z_j)+\beta}\setminus \stackrel{o}{B}_j,\cdot_j)\to \pi_{n-1}(F_{f(z_j)+\beta},\cdot_j)\to \pi_{n-1}(F_\alpha,\cdot)$\\
or of the variation mappings\\ 
$\pi_{n-1}(F_{t_j+\beta},\partial F_{t_j+\beta},\cdot_j)\to \pi_{n-1}(F_{t_j+\beta},\cdot_j)\to \pi_{n-1}(X_r\cap f^{-1}(\Gamma_j),\cdot_j)\to \pi_{n-1}(F_\alpha,\cdot)$\\
both being defined in Theorem 2.1.2.\\

{\bf Proof:} This follows from Lemma 2.2.1 and Theorem 2.1.2, with $f-t_j$ instead of $f$.\\

In the case $n=2$ we can proceed in a similar way as in Theorem 2.0.5:\\ 

{\bf Theorem 2.2.3:} Suppose $n=2$ and $F_\alpha$ connected.\\
a) Let $a_{j1},\ldots,a_{jr_j}\in\pi_1(F_{f(z_j)+\beta},F_{f(z_j)+\beta}\setminus \stackrel{o}{B}_j,\cdot_j)$ be the images of generators of $\pi_1(F_{f(z_j)+\beta},\cdot_j)$. Let $z_{j0}=\cdot_j,z_{j1},\ldots,z_{js_j}$ be representatives of the connected components of $F_{f(z_j)+\beta}\setminus \stackrel{o}{B}_j$, and let\\ 
$b_{jl}\in \pi_1(F_{f(z_j)+\beta},F_{f(z_j)+\beta}\setminus \stackrel{o}{B}_j,\cdot_j)$ be represented by a path in $F_{f(z_j)+\beta}$ from $z_j$ to $z_0=\cdot_j$. Then the subgroup generated by $Im\,Var^j_1(f)$ is the normal subgroup of $\pi_1(F_\alpha,\cdot)$ generated by the images of 
$$a_{j1},\ldots,a_{jr_j},b_{j1},\ldots,b_{js_j}$$ under $Var_1^j(f):\pi_1(F_{f(z_j)+\beta},F_{f(z_j)+\beta}\setminus \stackrel{o}{B}_j,\cdot_j)\to \pi_1(F_{f(z_j)+\beta},\cdot_j)\to \pi_1(F_\alpha,\cdot)$.\\
b) Let $a_1,\ldots,a_r\in\pi_1(F_\alpha,\partial F_\alpha,\cdot)$ be the images of generators of $\pi_1(F_\alpha,\cdot)$. Let $z_0=\cdot,z_1,\ldots,z_s$ be representatives of the connected components of $\partial F_\alpha$, and
let $b_l\in \pi_1(F_\alpha,\partial F_\alpha)$ be represented by a path from $z_l$ to $z_0=\cdot$. Then the subgroup generated by $Im\,Var^j_1(f)$ is the normal subgroup of $\pi_1(F_\alpha,\cdot)$ generated by the images of $a_1,\ldots,a_r,b_1,\ldots,b_s$ under $Var^j_1(f)$, where this time\\
$Var^j_1(f):\pi_1(F_\alpha,\partial F_\alpha,\cdot)\to \pi_1(F_{t_j+\beta},\partial F_{t_j+\beta},\cdot_j)\to \pi_1(F_{t_j+\beta},\cdot_j)\to \pi_1(F_\alpha,\cdot)$.\\

This is proved as Theorem 2.0.5, it facilitates computation in the case $n=2$.\\

As in the case of homology, the results of this section will be applied to the nearby fibre (we renounce here to the global case for homotopy groups).\\

{\bf 2.3. Application: hyperplane section of the local nearby fibre}\\

Now we turn to the set-up of section 1.3. Assume that $rhd(X)=n\ge 3$ and that $X$ has no strata of dimension $1$. From Lemma 1.3.1 and Theorem 2.2.2 we obtain:\\

{\bf Theorem 2.3.1:} The mapping $\pi_{n-2}(F\cap\{g=\beta\},\cdot)\to \pi_{n-2}(F,\cdot)$ is an epimorphism whose kernel is generated by the images of the variation mappings\\
$\pi_{n-2}(F\cap\{g=g(z_j)+\beta'\},F_\cap\{g=g(z_j)+\beta'\}\setminus \stackrel{o}{B}_j,\cdot_j)$\\
$\to \pi_{n-2}(F\cap\{g=g(z_j)+\beta'\},\cdot_j)\to \pi_{n-2}(F\cap\{g=\beta\},\cdot)$\\
or of the variation mappings\\ 
$\pi_{n-2}(F\cap\{g=t_j+\beta'\},\partial F\cap\{g=t_j+\beta'\},\cdot_j)\to \pi_{n-2}(F\cap\{g=t_j+\beta'\},\cdot_j)$\\
$\to \pi_{n-2}(F\cap\{g=\beta\},\cdot)$.\\

{\bf 2.4. Consequences for local homotopy groups}\\

Let $X_\varepsilon$ be defined as in section 1.0. We want to discuss consequences for the local homotopy groups $\pi_k(X_\varepsilon,\cdot)$, $0<k\le n-1$, hence $n\ge 2$. \\
We assume that $({\cal C}_n)$ holds, cf. subsection 2.0.\\

{\bf I.} Assume $n\ge 3$.\\

a) For $0<k<n-2$, hence $n\ge 4$, we have
$$\pi_k(X_\varepsilon,\cdot)\simeq \pi_k(X_\varepsilon\cap H,\cdot)$$
by the (weak) local Lefschetz-Zariski theorem ([HL2],[HL3]).\\

b) Now look at the case $0<k=n-2$, hence $n\ge 3$.\\
By the strong local Lefschetz-Zariski theorem (Theorem 1.0.1), we have\\
$\pi_{n-2}(X_\varepsilon,\cdot)$ $\simeq \pi_{n-2}(F,\cdot)$.\\
Then, $\pi_{n-2}(F,\cdot)$ can be calculated from $\pi_{n-2}(F\cap\{g=\beta\},\cdot)$ using  a general linear function $g$ and variation mappings, cf. Theorem 2.3.1.\\

c) Finally look at the case $0<k=n-1$, hence $n\ge 2$.\\
Here $\pi_{n-1}(X_\varepsilon,\cdot)$ can be calculated from $\pi_{n-1}(F,\cdot)$, using Theorem 2.0.3.\\

{\bf II.} In particular let us look at the fundamental group.\\

a) For $n\ge 4$ we have $\pi_1(X_\varepsilon,\cdot)\simeq \pi_1(X_\varepsilon\cap H,\cdot)$, by local Lefschetz-Zariski.\\

b) For $n=3$, $\pi_1(X_\varepsilon,\cdot)\simeq \pi_1(F,\cdot)$, and $\pi_1(F,\cdot)$ can be calculated from $\pi_1(F\cap\{g=\beta\},\cdot)$ by Theorem 2.2.3 and Theorem 2.3.1,  using variation mappings.\\
Note that $F\cap\{g=\beta\}$ can be assumed to be a Riemann surface, so\\
$\pi_1(F\cap\{g=\beta\},\cdot)$ is free, and we obtain a presentation of the original group.\\

c) For $n=2$, $\pi_1(X_\varepsilon,\cdot)$ can be computed from $\pi_1(F,\cdot)$ using a variation mapping (Theorem 2.0.5). Note that $F$ can be assumed to be a Riemann surface, so $\pi_1(F,\cdot)$ is free, and we obtain a presentation of the group in question.\\

{\bf References:}\\

[C] D. Ch\'eniot: Vanishing cycles in a pencil of hyperplane sections of a non-singular quasi-projective variety. Proc. London Math. Soc. (3) {\bf 72}, no. 3, 515-544 (1996).\\

[CE] D. Ch\'eniot, C. Eyral: Homotopical variations and high-dimensional Zariski-van Kampen theorems. Trans. Amer. Math. Soc. {\bf 358}, no. 1, 1-10 (2006).\\

[D] P. Deligne: Le groupe fondamental du compl\'ement 
d'une courbe plane n'ayant que des points  
doubles ordinaires est ab\'elien [d'apr\`es W. Fulton]. 
S\'eminaire Bourbaki {\bf 22}, 1-10 (1979-1980). \\

[E] C. Eyral: Homotopy theory on singular complex algebraic varieties. IMPAN Lecture Notes {\bf 5}. Instytut Matematyczny, Polska Akad. Nauk (IMPAN), 2022.\\

[H1] H. A. Hamm: Lokale topologische Eigenschaften komplexer R\"aume. Math. Ann. {\bf 191}, 235-252 (1972).\\ 

[H2] H. A. Hamm: A local Zariski-van Kampen theorem. J. of Sing. {\bf 30}, 216-241 (2026).\\

[HL1] H. A. Hamm, L\^e D. T.: Lefschetz Theorems on quasi-projective varieties. Bull. Soc. Math. France {\bf 113}, 123-142 (1985).\\

[HL2] H. A. Hamm, L\^e D. T.: Local generalizations of Lefschetz-Zariski theorems. J. reine angew. Math. {\bf 389}, 157-189 (1988).\\

[HL3] H. A. Hamm, L\^e D. T.: Rectified homotopical depth and Grothendieck conjectures. In: The Grothendieck Festschrift, Vol. II, pp. 311-351, Progr. Math. {\bf 87}, Birkh\"auser Boston, Boston, MA, 1990. \\

[L] S. Lefschetz: L'analysis situs et la g\'eom\'etrie alg\'ebrique. Reprint of the 1924 original. Gauthier-Villars, Paris, 1950.\\

[Sp] E. H. Spanier: Algebraic topology. Corrected reprint of the 1966 original, Springer-Verlag, New York, 1995.\\

[tD] T. tom Dieck: Algebraic topology. EMS Textbooks in Mathematics, 2008.\\

[Ti] M. Tib\v{a}r: Polynomials and vanishing cycles. Cambridge Tracts in Mathematics {\bf 170}, Cambridge Univ. Press, Cambridge, 2007.\\

[Z] O. Zariski: A theorem on the Poincar\'e group of an algebraic hypersurface. Ann. of Math. (2) {\bf 38}, 131-141 (1937).\\

Helmut A. Hamm\\
Mathematisches Institut\\
Universit\"at M\"unster, Germany\\
e-mail: hamm@uni-muenster.de

\end{document}